\documentclass[12pt]{article}

\usepackage{amsthm,amsmath,amssymb}

\usepackage{graphicx}
\usepackage[T1]{fontenc}

\usepackage[colorlinks=true,citecolor=black,linkcolor=black,urlcolor=blue]{hyperref}


\newcounter{parentnumber}
\theoremstyle{plain}
\newtheorem{theorem}{Theorem}

\theoremstyle{lemma}
\newtheorem{lemma}{Lemma}

\theoremstyle{proposition}
\newtheorem{proposition}{Proposition}

\theoremstyle{corollary}
\newtheorem{corollary}{Corollary}

\theoremstyle{definition}
\newtheorem{definition}{Definition}

\theoremstyle{remark}
\newtheorem{remark}[theorem]{Remark}

\title{\bf Reconstruction of graphs via asymmetry \\  \small{}  }

\author{Ameneh Farhadian \thanks{The contents of this paper are taken from the author's Ph.D. Thesis, Department of Mathematical Sciences, Sharif University of Technology, Supervised by Professor E. S. Mahmoodian} \\
\small Department of Mathematical Sciences\\[-0.8ex]
\small Sharif University of Technology\\[-0.8ex] 
\small  P. O. Box 11155-9415, Tehran, I. R. Iran\\
\small\tt a\_farhadian@mehr.sharif.ir\\ }

\date{ \small Mathematics Subject Classifications:  05C60}

\begin{document}

\maketitle


\begin{abstract}

Any graph which is not vertex transitive has a proper induced subgraph which is  unique due to its structure or the way of its connection to  the rest of the graph. We have called such subgraph as an \textit{anchor}. Using an  anchor which, in fact, is representative of a graph asymmetry, the reconstruction of that graph reduces to a smaller form of the reconstruction.
Therefore, to show that a graph is reconstructible, it is sufficient to find a suitable anchor that brings us to a solved form of the problem.
Let $G$ be an arbitrary graph with $n$ vertices which is not vertex-transitive. Usually,  graph $G$  has either a ($n-2$)-vertex anchor $H$ or an orbit $O$ with at least three vertices such that $G \backslash O$ is an anchor (or connective anchor) and it enables us to show that graph $G$ is reconstructible. For instance, this fact is enough to show that trees are reconstructible.

%
%

\end{abstract}

\section {Introduction}
The Reconstruction Conjecture is an interesting problem which has remained open for more than $70$ years. It states that all graphs on at least three vertices are determined up to isomorphism by their deck \cite{k, MR0120127}. The deck of a graph $G$ is the multiset of graphs that is obtained from deleting one vertex in every possible way from the graph $G$. The members of a deck are referred to as \textit{cards}. A class of graphs is \textit{reconstructible}, if every member of the class is reconstructible.\\
First time, this conjecture was proposed by Ulam \cite{ MR0120127}. Kelly \cite{MR0087949}, in his Ph.D thesis, showed that regular graphs, Eulerian graphs, disconnected graphs and trees are reconstructible. Bondy and Hemminger \cite{bondy1977graph} have proposed another proof for trees by employing the counting theorem.  Graphs in which no two cycle have a common edge \cite{manvel1969reconstruction}, graphs in which all  cycles pass through a common vertex \cite{manvel1978nearly}, Outer planner graphs \cite{giles1974reconstruction}, separable graphs without end vertices\cite{bondy1969ulam}, maximal planar graphs \cite{MR615314}, critical blocks and graphs with some specific degrees sequence \cite{Nash} are some well known classes of reconstructible graphs. In this paper, unique subgraph is employed to find new families of reconstructible graphs. Moreover, a general approach to graph reconstruction is proposed.\\
A unique subgraph which, in fact, is representative of a graph asymmetry, is a known key concept for the graph reconstruction problem. Bollob\'{a}s \cite{MR1037416} has employed graphs in which all ($n-2$) and ($n-3$)-vertex subgraphs are unique to show that almost every graph is reconstructible by three cards. Muller \cite{muller1976probabilistic} has shown almost every graph is reconstructible using graphs whose all ${n}/ {2}$-vertex subgraphs are unique. Unique subgraphs also have been used by Chinn \cite{chinn1971graph} and Zhu \cite{MR1470804} to introduce some families of reconstructible graphs. Ramachandran \cite{MR716448}, also, has employed the idea of unique subgraph for the digraph $N$-reconstruction.

In this paper, a general framework for the graph reconstruction problem is proposed. We show that an orbit $O$ of a graph $G$ which makes $G \backslash O$ to be an anchor or two vertices which makes $G \backslash \lbrace v,w\rbrace $ to be an anchor with the conditions which will be mentioned, is sufficient to show that $G$ is reconstructible.  This simple statement is sufficient to show many families of graphs are reconstructible. For instance, this fact is enough to show trees and small graphs are reconstructible. \\
After definition of anchor and shadow graph in the next section,  the reconstruction of graphs using anchor is introduced in the third section. Then, the concept of unique subgraph extension and maximal unique subgraph are employed to draw a general framework for the graph reconstruction problem in Section \ref{framework}. The reconstruction of graphs with ($n-2$)-vertex anchor is discussed in Section \ref{n-2}. Finally, as an application of the suggested framework, it is show that trees are reconstructible. A brief digest of the paper is given in the last section. 

\section{Definitions and Notations}
In this section, graph anchor and  shadow graph are defined and their relation is introduced.
In this paper, any subgraph is a vertex induced subgraph, otherwise, it is mentioned. The neighbors of any vertex $v$ in a graph is denoted by $N(v)$.
A graph $ G $ is called \textit{asymmetric}, if $\rm{Aut }(G)=I$. Two vertices $u$ and $v$ of a graph $G$ are called \textit{similar}, if there is an automorphism of $ G$ which maps $u$ into $v$. Dissimilar vertices whose removal leaves isomorphic subgraphs are called
\textit{pseudo-similar}\cite{lauri1997pseudosimilarity, JGT:JGT3190050207}. Similarity is, obviously, an equivalence relation. Thus, the similar vertices are in classes which are called \textit{orbits}. If $H$ is an subgraph of graph $G$, $G \backslash V(H)$ is induced subgraph on $V(G)-V(H)$.

\begin{figure}[ht]
\centerline{\includegraphics[width=9cm]{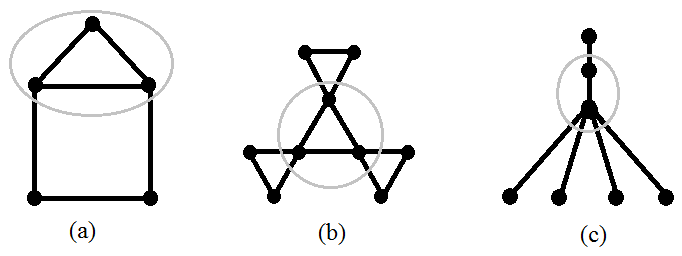}}
\caption{\label{connection}\small a)Anchor b,c)connective anchor}
\end{figure}

\begin{definition}
A proper induced subgraph $H$ of a graph $G$ is an \textit{anchor}, if it occurs exactly once in $G$.
 A subgraph $H$ of a graph $G$ which is not necessarily unique, but is distinct due to its connection to $G\backslash V(H)$ is a \textit{connective anchor}.
\end{definition}
The anchors of some graphs are shown in Fig. \ref{connection}. To have more intuition about anchors of a graph, please see the next section.

Unique subgraph, also, has been defined by Entringer and Erd\H{o}s \cite{entringer1972number} and used by Harary and Schewenk \cite{harary1973number}. They have used the concept of unique subgraph as the spanning subgraphs which are unique. Here, in contrast, we deal with proper vertex induced subgraph.   

An anchor is a unique subgraph and, therefore, is distinct in any card containing it. Therefore, an anchor in a graph, like a real anchor which fixes a boat, fixes a part of some cards and makes it possible to compare them.
\begin{definition}
A graph is \textit{balanced}, if  it does not have any anchor. A graph $G$ is quasi-balanced, if for any anchor $H$ of $G$, $\vert V(H) \vert =\vert V(G) \vert -1$.
\end{definition}A balanced graph may have a connective anchor. We will see in the following, a balanced graph that does not have any connective anchor is vertex-transitive.

Now, we define shadow graph. Using reconstructible families of shadow graphs, new families of reconstructible graphs will be determined which have not been recognized, previously.
\begin{definition}
Let $H(V,E)$ be an arbitrary graph and  $G(V_G,E_G)$  be a graph in which any vertex $v \in V_G$ is an arbitrary subset of $V_H$.  $E_G$  does not depend on $H$. We call $G$ is a \textit{shadow graph} on the background graph $H$.\\
 Two shadow graphs $G$ and $G'$ are \textit{isomorphic}, if their background graphs, that is $H$ and $H'$, are isomorphic and there is an isomorphism mapping $f: V(H) \to V(H')$ such that $f(V(G))=V(G')$ and any two vertices $u$ and $v$ of $G $ are adjacent in $G$ if and only if $f(u)$ and $f(v)$ are adjacent in $G'$.\\
To avoid confusion between ordinary graphs and shadow graphs, we use under-line for shadow graphs and its vertices and edges.
 An \textit{automorphism} group of a shadow graph $\underline{G}$ is an isomorphism mapping of $\underline{G}$ to itself. Two vertices of a shadow graph are \textit{similar}, if there is a shadow graph automorphism which maps one of them to another. A \textit{vertex-transitive} shadow graph is a shadow graph whose any two vertices are similar. If some vertices of a shadow graph are removed, we have a \textit{sub shadow graph}.  Two sub shadow graphs are isomorphic, if there is an automorphism of $H$ which maps vertices and preserve the adjacency.  A sub shadow graph is an \textit{anchor}, if there is not any sub shadow graph isomorphic to it. Please note that two vertices of a shadow graph may be non-isomorphic.  A vertex $\underline{v}$ of a shadow graph is \textit{fixed}, if $\theta (\underline{v})=\underline{v}$ for any $\theta \in Aut(H)$.
\end{definition}
 \begin{figure}[ht]\label{shadow}
\centerline{\includegraphics[width=8cm]{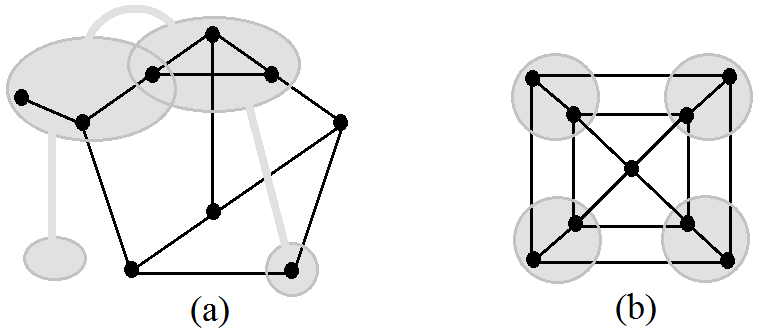}}
\caption{\label{Anchor_sup}\small (a) A shadow graph with four vertices and three edges (b) A vertex-transitive shadow graph with four vertices}
\end{figure}
In Fig. \ref{shadow}, two shadow graphs are shown. The right one is vertex-transitive. The background graphs are illustrated by black color and shadow graphs by gray color. Please note that the automorphism group of the vertex-transitive shadow graph (b) is $D_4$.
\section{ Anchor and Graph Reconstruction }
Here, we will be more familiar with anchors of a graph. Then, we use the anchor of a graph to reduce the reconstruction of that graph to the reconstruction of a smaller shadow graph in order.\\
 The concept of shadow graph provides a tool to compare the connection of two subgraphs. Let $H$ be a subgraph of a graph $G$ with two copies, say $H_1$ and $H_2$. To be able to compare the connection of $H_1$   and $H_2$ to the rest of the graph, we can compare two shadow graphs $G_{/H_1}$ and $G_{/H_2}$. 
\begin{lemma}
 If $O$ is an arbitrary orbit of a non-vertex transitive graph $G$, then $G \backslash O$ is unique due to either its structure or the way of its connection to $O$. 
\end{lemma}
\textit{Proof:} Let $H$ be $G \backslash O$. If $H$ is a unique subgraph, then it is finished. Thus, suppose that there is another copy of $H$ in  graph $G$, say $H'$. The shadow graph $G_{/H}$ is not isomorphic to $G_{/H'}$. Because, if two shadow graphs are isomorphic, then the vertices in $G \backslash V(H)$, i.e. $O$, are similar to vertices of $G \backslash V(H')$, while it contradicts to the assumption that  $O$ is an orbit. Because, an orbit includes all similar vertices of a graph.$\diamond$

\begin{corollary}
Any non-vertex-transitive graph has either an anchor or a connective anchor.
\end{corollary}

\begin{remark}
If a graph $G$ has a unique subgraph, then it can be determined from the deck.
Because, according to Kelly's Lemma, the number of
occurrence of each induced subgraph $H$ in $G$ on at most
$n-1$ vertices is reconstructed from the deck of $G$. Thus, any anchor is distinguishable in any card containing it. But, a connective anchor may be not recognizable in the deck. Because, if a vertex out of the anchor is deleted, we have a part of the connection of $H$ to the rest of the graph.
\end{remark}

\begin{figure}[ht]
\centerline{\includegraphics[width=9cm]{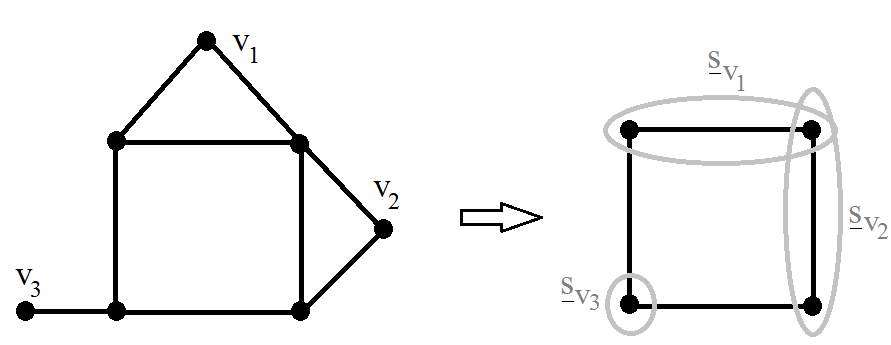}}
\caption{\label{Anchor_sup}\small The left graph with anchor $C_4$ is converted into the shadow graph in the right side.}
\end{figure}

 Now, we explain the relation of ordinary graphs and shadow graphs. 
For any subgraph $H$ of a graph $G$, a shadow graph can be dedicated in which $H$ is a background graph. Let $H$ be an anchorof a graph $G$. We take the anchor $H$ as the background graph and define  $\underline{s_v}= N(v) \cap V(H)$  as a vertex of a shadow graph for any $v \in V(G)-V(H)$. We make two vertices $\underline{s_v}$ and $\underline{s_w}$ adjacent in the shadow graph, if they are adjacent in $G$. Therefore, graph $G$ is converted to a shadow graph on the background graph $H$. We denote this shadow graph by $\underline{G}_{/H}$.  In Fig.\ref{Anchor_sup}, it is shown that how a graph turns into a shadow graph. The vertices of the shadow graph are shown by gray color. This conversion is reversible.

The set of cards containing the anchor forms the deck of the shadow graph established on the anchor. If the shadow graph is reconstructible from deck, then the graph is reconstructible.
\begin{proposition}\label{base}
Let $H$ be an anchor of a graph $G$ with at most ($n -3$) vertics. Graph $G$ is reconstructible, if the shadow graph $\underline{G}_{/H}$ is reconstructible.
\end{proposition}
\textit{Proof}: The anchor $H$ can be determined from the deck. Because, the definition of an anchor is that it
appears exactly once as an induced subgraph of $G$,
and Kelly's Lemma tells us that the number of
occurrence of each induced subgraph $H$ in $G$ on at most
$n-1$ vertices can be reconstructed from the deck of $G$. Thus, any anchor is distinguishable in any card containing it.
 The cards containing the anchor establish the deck of shadow graph $\underline{G}_{/H}$. If this shadow graph is reconstructible, the neighbors of each vertex of $G \backslash V(H)$ is determined. Thus, $G$ is reconstructible. $\diamond$\\
Please note that a shadow graph is not necessarily reconstructible. But, reconstructible shadow graphs provide significant families of  reconstructible graphs.
\begin{definition}
A connective  anchor is distinguishable, if it can be distinguished in the deck.
\end{definition}
 A distinguishable connective anchor, like an anchor, is distinct in any card containing it and plays the same role of an anchor for the reconstruction of a graph from its deck.

 It may be claimed that since  just cards containing the anchor are used for graph reconstruction, this method can not be efficient. But, this claim is not exact. Because, to find the anchors of a graph, we should consider all subgraphs from all cards to find the unique ones. In addition, the anchor extension idea which is introduced in the next section, leads us to employ all cards of the deck for the graph reconstruction. Therefore, in this method we implicitly employ all cards for reconstruction. We will see in the next section, how the anchor concept makes a breakthrough in the graph reconstruction problem. 
\section{A Framework for the Graph Reconstruction Problem }\label{framework}
In this section, the concept of anchor extension and maximal anchor are introduced which draw a general framework for the graph reconstruction problem.  It states that the  reconstruction conjecture is true, if we can prove it for a) balanced and qusi-balanced graphs and shadow graphs with at least three vertices and b) graphs with a $(n-2)$-vertex anchor.

\begin{definition}
An anchor is maximal, if it can not be extended to a larger anchor.
\end{definition}
\begin{lemma}\label{maximal}
If $H$ is a maximal anchor of a graph $G$, then the shadow graph $\underline{G}_{/H}$  is balanced.
\end{lemma}
\textit{Proof}: Suppose that the shadow graph $\underline{G}_{/H}$ has an anchor and let the induced sub shadow graph on $A \subset V( \underline{G }_{/H}) $ be its anchor. We show the induced subgraph on $V(H) \cup A$ is an anchor of graph $G$, contradicting to the maximality of anchor $H$. Suppose that the induced subgraph on $V(H) \cup A$ is not an anchor of $G$. Thus, there is $B \subset V(G)-V(H)$ that the induced subgraph on $V(H) \cup A$ and $V(H) \cup B$ are isomorphic by mapping $\phi : V(H) \cup A \to V(H) \cup B$. Since, $H$ is an anchor and unique, $\phi(H)=H$ and $\phi(A)=B$. Thus, the restriction of $\phi$ to $V(H)$ is an automorphism of $H$ and the induced shadow graph on $A$ is isomorphic to the induced shadow graph on $B$, contradicting to $A$ to be an anchor of the shadow graph $\underline{G }_{/H}$.$\diamond$

The above lemma has useful results. It, in  fact, provides an algorithm to investigate a graph for reconstruction. According to this lemma, any anchor of a graph can be extended until we do not have any anchor out of the anchor. In other words, this lemma reduces the reconstruction of a graph to an balanced structure or a graph with a sufficiently large anchor, i.e. an $(n-2)$-vertex anchor. Therefore, it provides a general framework for the reconstruction problem which is drawn in the following theorem.

\begin{theorem}\label{main}
The conjecture of graph reconstruction is true, if we prove the following two families of graphs are reconstructible:\\
a)  graphs and shadow graphs which are balanced or quasi-balanced with at least three vertices and\\
b) $n$-vertex graphs with ($n-2$)-vertex anchor.
\end{theorem}
\textit{Proof}: Let $G$ be an arbitrary graph with at least three vertices. We want to show $G$ is reconstructible, provided that the families of (a) and (b) are reconstructible. If $G$ is balanced or quasi-balanced, then it is reconstructible due to (a). Thus, suppose that $G$ has an anchor which by extending it, we have three cases: \\
Case \textit{I}: We reach to a maximal anchor $H_m$ with at most ($n-3$) vertices. Thus, the shadow graph $\underline{G}_{/ H_m}$ is balanced due to Lemma \ref{maximal} and is reconstructible due to assumption (a). Therefore, $G$ is reconstructible due to proposition \ref{base}.\\
Case  \textit{II}: We reach to an ($n-2$)-vertex anchor which is reconstructible by the assumption (b). \\
Case  \textit{III}: We reach to an ($n-1$)-vertex anchor $H$. We have supposed that $G$ is not quasi-balanced, thus we have  reached to $H$ by the extension of an anchor, say $H'$. If $\underline{G}_{/H'}$ has an anchor of order $k$ where $k<\vert V( \underline{G}_{/H'} )\vert-1 $, then we can add this anchor to $H'$ to have a larger anchor and reach to the above cases of \textit{I} or \textit{II}. But, if $\underline{G}_{/H'}$ does not have any such anchor, it is a quasi-balanced shadow graph which is reconstructible by the assumption (a). Thus, $G$ is reconstructible.
$ \diamond $ 

The above theorem, in fact, proposes an efficient approach to investigate any graph for reconstruction from its deck. According to this theorem, if we show that the balanced graphs and shadow graphs and, also, graphs with an ($n-2$)-vertex anchor are reconstructible, then we have solved the reconstruction problem. THese two families of graphs are investigated in the two next sections.

\section{Balanced and Quasi-Balanced Graphs}
The majority of balanced graphs and shadow graphs are vertex-transitive. We know that vertex-transitive graphs are reconstructible due to regularity. Here, we show vertex-transitive shadow graphs are, also, reconstructible. Then we deal with quasi-balanced graphs.

\begin{lemma}\label{s1}
If shadow graph $\underline{G}$ with at least three vertices is vertex-transitive, then $\underline{G}$ is reconstructible.
\end{lemma}
\textit{Proof}: Shadow-graph $\underline{G}$ is vertex-transitive. Thus, all cards of the deck are the same. The neighbors of $\underline{v}$ is recognizable by regularity. It is sufficient to find the location of the vertex $\underline{v}$ in the background graph $H$, i.e. subset $\underline{v}$ of the set $V(H)$.
Since $\vert V(O) \vert > 3$, there are at least two different vertices, say $\underline{w}$  and $\underline{w'}$, in $V(\underline{G})$. Since they belong to the same orbit, there exists $\alpha \in Aut(H)$  such that $\alpha(\underline{w})=\underline{w'}$. 
We assign a digraph to the shadow graph.  For any $\underline{w}, \underline{w'} \in V(\underline{G})$ such that $\underline{w'}= \alpha (\underline{w})$, we draw an arc from $\underline{w}$ to $\underline{w'}$. The resulted digraph is vertex transitive. Because, all cards are the same. Thus, it is reconstructible due to regularity. Thus, in card $\underline{G} \backslash \underline{v}$, $\underline{v}$ is $\alpha(\underline{w})$ for a vertex $\underline{w} \in V(\underline{G} \backslash \underline{v})$.$\diamond$\\
\begin{figure}[ht]
\centerline{\includegraphics[width=6.7cm]{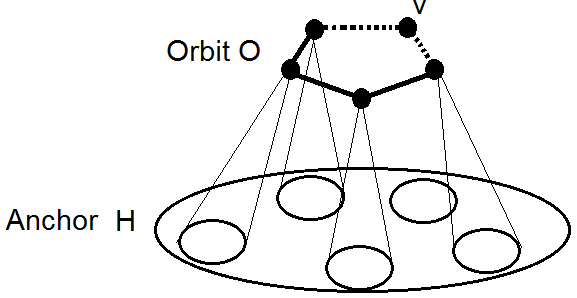}}
\caption{\label{Anchor_example}\small Graph $G$ with orbit $O$ and anchor $G\backslash O$ is reconstructible.}
\end{figure}
The following theorem is a corollary of the above Lemma.
\begin{theorem}\label{thm:twoB}
Let $O$ be an orbit of a graph $G$ with at least three vertices. If $G \backslash O$ is an anchor, then $G$ is reconstructible.
\end{theorem}
\textit{Proof}: It is sufficient to take the established shadow graph on the anchor $G \backslash O$. The result follows from Lemma \ref{s1} and preposition \ref{base}.$\diamond$\\
In the above theorem, $G \backslash O$ can be a connective anchor which is distinguishable in the deck.

\begin{figure}[ht]
\centerline{\includegraphics[width=14cm]{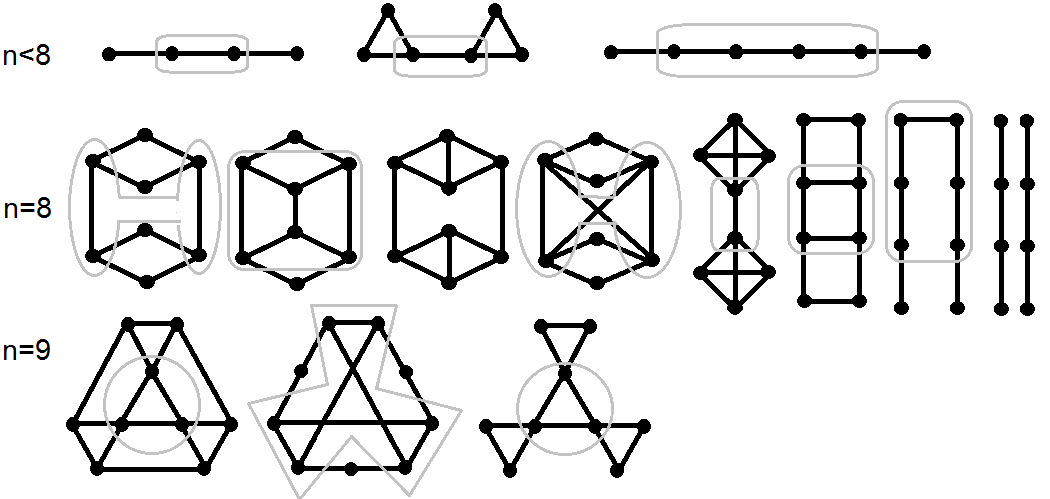}}
\caption{\label{anchor free}\small All balanced graphs which are not vertex transitive with at most 9 vertices. The connective anchors that prove they are reconstructible are shown.}
\end{figure}
Although the majority of balanced graphs are vertex-transitive, but there are special cases which are not vertex-transitive. All balanced graphs with at most 9 vertices, which are not vertex transitive, are shown in Figure \ref{anchor free}. Since both of a graph and its complement are balanced or not, a graph or its complement is considered in the figure. We see that for all of them, the deletion of an orbit  leaves a connective anchor which is distinguishable in the  deck. The connective anchor of them are shown by a gray curve. Therefore, all of them are reconstructible.\\
\begin{figure}[ht]
\centerline{\includegraphics[width=12cm]{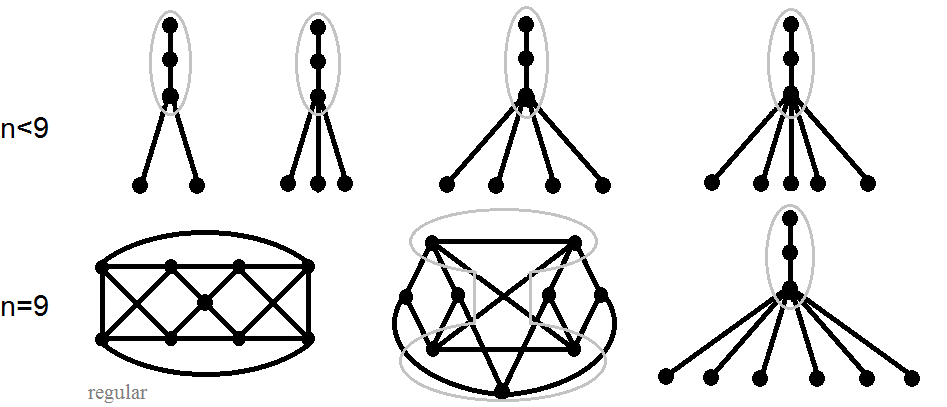}}
\caption{\label{n-1}\small  All quasi-balanced graphs with at most 9 vertices  which there is not any isolated vertex in them or their complement. }
\end{figure}
 The majority of quasi-balanced graphs are disjoint union of an isolated vertex with a balanced graph or its complement. Thus, these quasi-balanced graphs are reconstructible due to be disconnected. The exceptional cases with at most 9 vertex which do not have such structure, are shown in Figure \ref{n-1}. Again, a graph or its complement are  shown. The distinguishable connective anchor of these graphs are, also, shown in the figure.

\section{Graphs with  ($n-2$)-vertex anchor}\label{n-2}
Now, we discuss about the reconstruction of graphs with an ($n-2$)-vertex anchor which was mentioned in the part(b) of Theorem \ref{main}.

Chinn \cite{chinn1971graph} has shown if there exists a vertex $v$ such that all ($n-2$)-vertex subgraphs of $G \backslash \lbrace v \rbrace$ are unique, graph $G$ is reconstructible. Zhu \cite{MR1470804} has improved this result by showing that at most three of ($n-2$)-vertex subgraphs of  $G \backslash \lbrace v \rbrace$  can be non-unique. Here, we show that one unique ($n-2$)-vertex subgraph which is asymmetric, is enough for $G$ to be reconstructible.

\begin{theorem}\label{fix}
Let $G$ be a graph with unique subgraph $H=G \backslash \lbrace v, u \rbrace$.  If $\alpha (N_{v,H})=N_{v,H}$ for any $ \alpha\in Aut(H)$, then $G$ is reconstructible.
\end{theorem}
\textit{Proof}: The anchor $H$ is distinct in the cards  $G \backslash \lbrace v \rbrace$ and  $G \backslash \lbrace u \rbrace$. Thus, we find the neighbors of $v$ in $H$ using the card  $G \backslash \lbrace u \rbrace$. Since, $\alpha (N_{v,H})=N_{v,H}$ for any $ \alpha\in Aut(H)$, there is just one way to add $v$ to $H$ in the card  $G \backslash \lbrace v \rbrace$. The existence of edge between $v$ and $u$ can be inferred from the number of edges which is reconstructible.$\diamond $

\begin{corollary}\label{asymn-2}
Any $n$-vertex graph with an asymmetric unique subgraph of order $n-2$ is reconstructible.
\end{corollary}
The above corollary satisfies for almost every graph. Because, in almost every graph all $(n-3)$-vertex subgraphs are mutually non-isomorphic \cite{korshunov1985main, muller1976probabilistic, MR1037416}. In such graphs, any ($n-2$)-vertex subgraph is unique and asymmetric. 

\begin{lemma}\label{Z2}
Let $O$ be an orbit of a graph $G$ with two vertices such that $G \backslash O$ is an anchor and Aut($G \backslash O$) $\in \lbrace I, Z_2, Z_3, D_3 \rbrace$, then graph $G$ is reconstructible.
\end{lemma}
\textit{Proof}: If Aut($G \backslash O$)=$I$, then $G$ is reconstructible due to the previous corollary. If Aut($G \backslash O$) $\in \lbrace Z_2, Z_3, D_3 \rbrace$,  there are two possibilities for the neighbors of $ V(O)$ in $G \backslash O$. Two vertices of $O$ have either the same neighbors in $G \backslash O$ or different neighbors. The degree sequence of vertices which is reconstructible by the deck, separates these two cases and chooses one of them.$\diamond$\\
In a graph $G$ with anchor $H=G \backslash \lbrace v,w \rbrace$, let $s_v$ and $s_w$ be the neighbors of $v$ and $w$ in $H$, respectively. $s_v$ and $s_w$ are recognizable from cards $G \backslash w$ and $G \backslash v$, respectively. But, $s_v$ and $s_w$ can arbitrary move by the action of $H$ automorphism group. It may make different possibilities for the state of placing both of them together and, consequently, two cards containing the anchor are not sufficient to reconstruct the graph.

\begin{figure}[ht]
\centerline{\includegraphics[width=10cm]{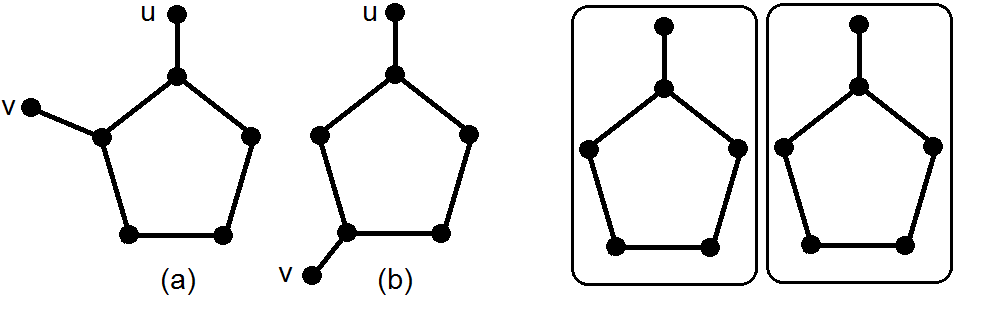}}
\caption{\label{relative}\small For both graphs (a) and (b), the vertices $\underline{s_{v}}$ and $\underline{s_{w}}$ are isomorphic in the shadow graph established on the anchor $C_4$. But, their relative position makes non-isomorphic graphs. }
\end{figure}

For example, in Fig. \ref{relative}, two non-isomorphic graphs a and b are shown. The subgraph $C_5$ is an anchor for both of them. The set of cards containing the anchor are the same and is shown in the right side. The set of cards containing the anchor are not sufficient to discriminate these two non-isomorphic graphs. Therefore, it is necessary to use other cards to find their relative  position. We do not know whether the other $n-2$ cards are, always, sufficient to determine the state of putting both of them together. The reconstruction conjecture claims that they are sufficient. 

The following theorem uses the distance of vertices $v$ and $w$ to identify the state of placing both of their neighbors together, when the automorphism group of the anchor is not trivial. For example, in  Fig. \ref{relative} the distance between $v$ and $w$ in two graphs discriminates these two graphs. Thus, these two graphs are reconstructible due to the next theorem.

\begin{lemma} \label{distance_anchor}
Let $G$ be a graph with anchor $H=G \backslash \lbrace v, u \rbrace$. If the distance of $v$ and $w$ within $H$ specifies the state of placing both of the neighbors of $ v$ and $w$ together in $H$, then graph $G$  is reconstructible.
\end{lemma}
\textit{Proof}: The neighbors of $w$ and $v$ in $H$  are distinguishable in cards $G\backslash v$ and $G\backslash w$, respectively, up to isomorphism. It is sufficient to know the position of them when they come together. According to hypothesis, the distance of $v$ and $w$ within $H$ clarifies the relative position of the neighbors of $v$ and $w$ in $H$. Thus, it is sufficient to show the distance of $v$ and $w$ within $H$ is reconstructible.\\
The number of subgraphs containing $v$ and the number of subgraphs containing $w$ can  be obtained from the cards $G\backslash w$ and $G\backslash v$, respectively. In addition, the subgraphs which include  none of them are, exactly, the subgraphs of $H$ and, thus, their numbers is reconstructible. Therefore, the number of subgraphs which include both $v$ and $w$ are reconstructible. The smallest path which includes $v$ and $w$ indicates the distance of $v$ and $w$ in $H$ when $v$ and $w$ are not adjacent in $G$. If $v$ and $w$ are adjacent in $G$, we consider the smallest cycle including $v$ and $w$.  Thus, the distance of $v$ and $w$ within $H$ is reconstructible. $\diamond$

 The following corollary shows that the above Lemma always satisfies for tree leaves.
\begin{corollary} \label{tree2}
Let $T$ be a tree with anchor $T_{anch}=T\backslash \lbrace v,w \rbrace$ where $v$ and $w$ are two leaves of $T$, then $T$ is reconstructible. 
\end{corollary}
\textit{Proof}:
We demonstrate that the distance of $v$ and $w$ in $T$ specifies the state of placing both of them together. Thus, $T$ is reconstructible due to the previous lemma.\\
For a leaf $x$ out of anchor $T_{anch}$, let $s_x$ be the neighbor of $x$ in $T_{anch}$.\\
Let $T_{anch}$ be the background graph. For any vertex $v$ out of anchor, $s_v$ is a vertex of a shadow graph. Thus, we have a shadow graph on $T_{anch}$.
We show that $\lbrace s_v, s_w \rbrace $ on $T_{anch}$ is  isomorphic to the $\lbrace s_v' ,s_w' \rbrace$ on $T_{anch}$ if and only if $d(v,w)=d(v',w')$. From left to right is trivial. Thus, we prove the converse.\\
We want to show if $d(v,w)=d(v',w')$, then there is a mapping in automorphism of $T_{anch}$ where maps $\lbrace s_v, s_w \rbrace $ to $\lbrace s_{v'}, s_{w'} \rbrace$. Since $s_v$ and $s_{v'}$ belongs to the same orbit, without loose of generality, we assume $v=v'$. Now, it is sufficient to show that there is a mapping in stabilizer of $s_v$ where maps $s_w$ to $s_{w'}$. 
Let $P : s_v = v_0, v_1,...,v_k = s_w$ be the $s_v$-$ s_w$ path, $Q$ the $s_v$ -$ s_{w'}$ path, and $t$ the largest integer ($0<t <k -1$) for which $v_t \in  V (P ) \cap V (Q) $. Necessarily, $d(s_w, v_t) = d(s_{w'}, v_t)$. Let $T_w$ and $T_{w'}$ be the  components containing $w$ and $w'$, respectively, in $T_{anch}\backslash v_t$. According to \cite{prins}, there is an automorphism of $T_{anch}$ that interchanges the components $T_w$ and $T_{w'}$, while other vertices are fixed. Since the vertex $s(v)$ does not belong to neither $T_w$ nor  $T_{w'}$, it is fixed under such mapping. Therefore, this mapping belongs to the stabilizer of $s(v)$. $\diamond$
\section{An application of the suggested framework}\label{app}

As an application of this framework, we will use it to show trees are reconstructible. The application of anchor for the reconstruction of the small graphs are given in Appendix I.

First time, Kelly \cite{k}, in his PhD thesis, proved that trees are reconstructible. His proof is relatively long. Bondy and Hemminger \cite{bondy1977graph} have proposed another proof for this fact by employing counting theorem. Here, a new proof is given for trees to be reconstructible.

\begin{corollary}\label{}
Trees are reconstructible. \cite{k}
\end{corollary}
\textit{Proof}:
 Let $T$ be a tree and $H $ be the induced subgraph on the vertices which are not tree leaves. We show $H$ is an anchor or a connective anchor which is distinguishable. If $H$ is a unique subgraph, then it is an anchor. Otherwise, $H$ is unique due to its connection, because the number of edges between $H$ and $T \backslash H$ is maximum possible value. If $T$ has $k$ leaves, in any card that a leaf is deleted from, we have $k-1$ or $k$ leaves.  If there are $k-1$ leaves, $H$  can be obtained by deletion of $k-1$ leaves and if there are $k$ leaves, deletion of any subset of $k-1$ leaves that results in a subgraph isomorphic to $H$ gives $H$. Because, if deletion of two leaves of a tree results in two isomorphic subgraph, then those two leaves are similar. Therefore, $H$ is distinguishable in any card containing it.\\
Any tree has at least two leaves. Thus, $T \backslash H$ has at least two vertices. If $T-H$ has two vertices or is an orbit with at least three vertices, then $T$ is reconstructible, respectively, due to Corollary \ref{tree2} and Theorem \ref{thm:twoB}. Otherwise, we enlarge the anchor  until we reach to two vertices out of the anchor or a maximal anchor with at most $n-2$ vertices. If we reach to a maximal anchor, the shadow graph established on the anchor is balanced due to Lemma \ref{maximal}. If only one vertex remains in $T\backslash H$, we come back one step, and  have a quasi-balanced shadow graph.\\
We show this balanced or quasi-balanced shadow graph is reconstructible and, consequently, $T$ is reconstructible. If the vertices of the shadow graph belong to two or more classes of isomorphism, then the induced sub shadow graph on one class make an anchor. Thus, all vertices of the shadow graph are isomorphic. Let $v_c$ be the center of $H$. We assign to each vertex $u \in V(H)$, the number of vertices in $V(T)-V(H)$  which their path to  $v_c$ pass through the vertex $u$.
In any depth, all non-zero values should be the same ( or at most there is just one vertex with value of 1 and other are the same when shadow graph is quasi-balanced). Because otherwise, vertices with maximum values enables us to have an anchor. Therefor, we can find the location of a vertex of shadow graph from disequilibrium of the assigned values in a card that a vertex of shadow graph is omitted from. In fact, in this case the shadow graph is vertex transitive or there is a vertex which its deletion reveals a vertex transitive shadow graph.$\diamond$
\section{Conclusion}
We saw that any graph which is not vertex-transitive has either an anchor or a connective anchor. We extend the anchor (or distinguishable connective anchor) as far as we  reach an orbit with at least three vertices or two vertices  out of the anchor. Then, we use the following result to show that the graph is reconstructible.

Let $G$ be a  graph which is not vertex-transitive with $n$ vertices. If $G$ has an orbit $O$ with at least three vertices  such that $G\backslash O$ is an anchor (connective anchor) or  has an anchor (connective anchor) with ($n-2$) vertices such that it satisfies in the condition of Theorem \ref{fix} or Lemma \ref{Z2} or Lemma \ref{distance_anchor}, then $G$ is reconstructible.

The above statement, in spite of being  simple, is sufficient to show that considerable families of graphs are reconstructible from their deck. For instance, using the above statement, we have shown that  trees and small graphs are reconstructible. Until now, we have not found any graph that does not satisfy in the condition of the above statement.
\section{Acknowledgment}
I would like to thank Prof. E. S. Mahmoodian and, also, Prof. A. Daneshgar for the time that they have dedicated to me.
 \bibliographystyle{plain}

\section*{Appendix I: Using anchor to show small  graphs are reconstructible}
\begin{figure}[ht]
\centerline{\includegraphics[width=15cm]{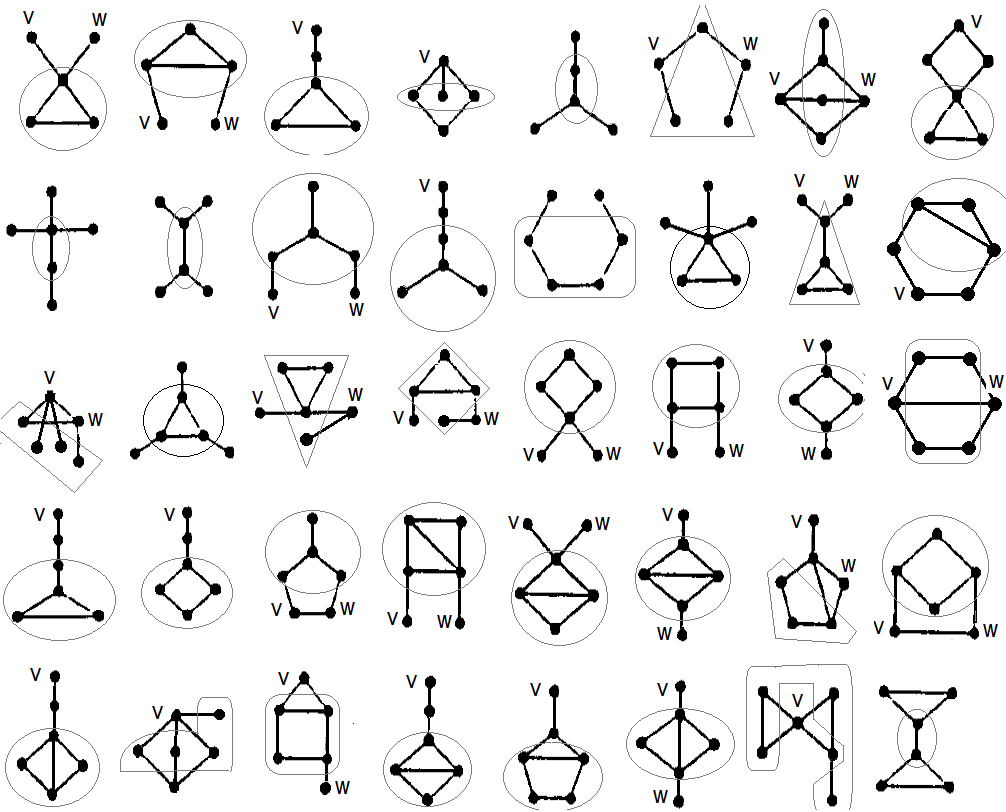}}
\caption{\label{small}\small All graphs with at most 6 vertices which are not disconnected or regular (also their complement) with the anchor that proves they are reconstructible.  }
\end{figure}

All graphs $G$ with at most 6 vertices such that $G$ and $G^c$ are not disconnected or regular, are shown in Fig. \ref{small}. Their anchor which  makes them  to be reconstructible are shown by a gray closed curve. 

\end{document}